\documentclass[11pt]{article}

\usepackage{amsmath,amsfonts,amscd, amssymb,theorem}

\newcommand{\doubleheaddownarrow}{\big\downarrow\kern-3.325mm\downarrow}

\newcommand\Z{\mathbb Z}
\newcommand\Q{\mathbb Q}
\newcommand\C{\mathbb C}

\newcommand\om{\omega}

\newcommand\Pf{\operatorname{Pf}}

\newcommand{\broken}{\dasharrow}
\newcommand{\blank}{\, \hbox{--} \,}

\newcommand{\id}{\operatorname{id}}

\newcommand\rank{\operatorname{rank}}

\newcommand\Ext{\operatorname{Ext}}

\newcommand\Hom{\operatorname{Hom}}
\newcommand\res{\operatorname{res}}

\newcommand{\skewentry}{ \kern -1cm -\mathrm{sym} \kern -1cm}
\newcommand\Spec{\operatorname{Spec}}

\renewcommand\P{\operatorname {\mathbb P}}
 \newcommand\Pp{\operatorname {\mathbb P}}

 \newtheorem{theorem}{Theorem}[section]
 \newtheorem{lemma}[theorem]{Lemma}

 {
 \theorembodyfont{\rmfamily}
 \newtheorem{defn}[theorem]{Definition}
 \newtheorem{exa}[theorem]{Example}

 \newtheorem{rem}[theorem]{Remark}

 }
 \newenvironment{pf}{\paragraph{Proof}}{\par\medskip}
 
 \newcommand{\qed}{\ifhmode\unskip\nobreak\fi\quad\ensuremath\square}
\newcommand{\QED}{\ifhmode\unskip\nobreak\fi\quad\ensuremath{\mathrm{QED}}}

\numberwithin{equation}{section}

\title{ Kustin--Miller unprojection with complexes }
 
\author{Stavros Papadakis}
\date {November 2001}

\begin{document}
\maketitle

\begin {abstract}
A main ingredient for Kustin--Miller unprojection, as developed in
\cite{PR}, is the module $\Hom_R(I, \om_R)$, where $R$ is
a local Gorenstein ring and $I$ a codimension one ideal with $R/I$ 
Gorenstein. We prove a method of calculating it in a relative 
setting using resolutions. We give three applications. In the 
first we generalise a result of \cite{CFHR}. The second and the third 
are about Tom and Jerry, two families of Gorenstein codimension
four rings with $9 \times 16$ resolutions.
\end{abstract}

\section {Introduction}

Kustin--Miller unprojection, as developed in \cite{KM} and \cite{PR}, is a
method 
for constructing a birationally equivalent Gorenstein scheme $Y$ to a
Gorenstein 
scheme $X$ (containing some codimension one subscheme $D$). A useful
variant
of this is a global version for a projectively Gorenstein scheme,
producing a 
birationally equivalent projective scheme. It has found many applications 
in algebraic geometry, 
for example in the birational  geometry of Fanos \cite{T}, \cite{CPR} and
\cite{CM}, 
in the construction of weighted complete  intersection K3s and  Fanos
\cite{Al},
and 
in the study of Mori flips \cite{BrR} (compare also \cite{Ki}). In all
cases,
the algebra of the problem was relatively simple and was handled with ad
hoc
methods.

In the present work we present systematic tools that calculate the
unprojection
in more 
complicated cases, by relating the conceptual approach based on the
adjunction
formula of \cite{PR} with the complex oriented  view of \cite{KM}. 

To demonstrate our methods, we study three important cases of
unprojection.
The first and simplest  (Section~\ref{sec!uciinsideci})
is the unprojection of a complete intersection inside
a complete intersection, where we generalise  a result of \cite{CFHR}.

The Tom and Jerry families of unprojection
(Section~\ref{sec!tomandjerry}), 
which were defined 
and named by Reid, play a more important role in algebraic geometry
because they
appear 
in hundreds of known examples; it is an open problem (\cite{Ki},
Problem~8.3)
whether
all the currently known  Gorenstein codimension~4 rings with $9 \times 16$
resolution
can be accommodated within these two structures. 

Our main results concerning Tom and Jerry are
Theorems~\ref{thm!localtomthm}
and~\ref{thm!localjerrythm}, where we calculate their unprojections using
multilinear and
homological algebra. An  interesting observation is  that they can be
considered 
as a kind of deformation of two standard models: the graded rings of the
Segre
embeddings of 
$\P^2 \times \P^2 \subset \P^8$ and of $\P^1 \times \P^1 \times \P^1
\subset
\P^7$ respectively;
these standard models play an important part in our arguments.  
For more open problems and further discussion of Tom and Jerry compare
\cite{Ki}~Section~8.

Finally, in Section~\ref{sec!fc} we state some general remarks about the
calculation of
unprojection. 

\section {Notation}   \label{sec!chpt1notat}

\paragraph {Gorenstein ideals}
In the present work, an ideal $I$ of a Gorenstein ring $S$ is called
Gorenstein
when the quotient ring $S/I$ is Gorenstein.

\paragraph {Pfaffians}
Assume $A=[a_{ij}]$ is a $k \times k$  skewsymmetric 
(i.e., $a_{ji}=-a_{ij}$ and $a_{ii}=0$) matrix with entries in 
a Noetherian ring $S$.  

For $k$ even we define a polynomial $\Pf(A)$  in $a_{ij}$  called  
{\em Pfaffian} of $A$  by induction on $k$. If $k=2$ we set
\[
     \Pf (
            \begin {pmatrix}  0 & a_{12} \\
                             -a_{12} & 0 
            \end {pmatrix} ) = a_{12}.
\] 
For even $k \geq 4$ we define
\[
     \Pf (A) = \sum_{j=2}^k (-1)^{j} a_{1j}\Pf (A_{1j}), 
\]  
where $A_{1j}$ is the skewsymmetric submatrix of $A$ obtained by deleting
the first and the $j$th row and column of $A$. An interesting property is
that 
\[
   (\Pf (A))^2  = \det A. 
\]  
 
Now assume that $k = 2l+1$ is odd. In the present work, by {\em Pfaffians}
of 
$A$ we mean the set 
\[
   \big\{ \Pf(A_1), \Pf (A_2), \dots , \Pf (A_k) \big\},   
\]
where for $1 \leq i \leq k$ we denote by $A_i$ the skewsymmetric submatrix
of $A$ obtained by deleting the $i$th row and column of $A$. 
Moreover, there is a complex  ${\bf L}$:
\begin{equation}   \label{eq!pfafcompl} 
  0 \to S \xrightarrow{\, B_3 \,} S^k \xrightarrow{\, B_2 \,} S^k 
      \xrightarrow{\, B_1 \,} S \to 0      
\end{equation}
associated  to $A$, with $B_2=A $, $B_1$ the $1 \times k$ matrix with
$i$th entry equal to $(-1)^{i+1}\Pf (A_i)$ and $B_3$  the transpose matrix
of $B_1$.  We will use the following theorem due to Eisenbud and Buchsbaum
\cite{BE}.

\begin {theorem} \label{thm!grade3compexact}
Let  $S$ be a Noetherian ring, $k=2l+1$ an odd integer and $A$ a 
skewsymmetric $k \times k$ matrix with entries in $S$. Denote by $I$
the ideal generated by the Pfaffians of $A$. Assume that $I \not= S$ and 
the grade of $I$ is three, the maximal possible. Then the complex ${\bf
L}$ 
is acyclic, in the sense that the complex 
\begin{equation*}   
  0 \to S \xrightarrow{\, B_3 \,} S^k \xrightarrow{\, B_2 \,} 
   S^k \xrightarrow{\, B_1 \,}  S \to S/I \to 0      
\end{equation*}
is exact. In addition, if $S$ is Gorenstein then the same is true for
$S/I$.
\end{theorem} 

For more details about Pfaffians  and a proof of the theorem  see 
e.g.~\cite{BE} or \cite{BH} Section~3.4.

\section {Calculation of $\Hom$ module }

Let $S$ be a Gorenstein local ring and $I \subset J$ perfect (therefore
Cohen--Macaulay)  
ideals of codimensions $r$  and $r+1$ respectively. Unless otherwise
indicated,
all $\Hom$ 
and $\Ext$ modules and  maps are over $S$. 

Since $S/I$ is Cohen--Macaulay, the adjunction formula (compare \cite{R1} 
p.~708 or \cite{AK} p.~6)  gives
\[
 \om_{S/J}=\Ext_{S/I}^1(S/J,\om_{S/I}). 
\]

We Hom (over the ring $S/I$) the exact sequence $0\to J/I\to S/I
\to S/J \to0$
with $\om_{S/I}$ to get the fundamental exact  sequence 
\begin {equation}   \label{eq!KMladjunct}
   0  \to \om_{S/I} \xrightarrow{\, a\,} \Hom_{S/I} (J/I, \om_{S/I}) 
             \xrightarrow{\, \res \,} \om_{S/J} \to 0,
\end{equation} 
with $a$ the natural map 
\[
        a(x)(l) = lx, \quad  \hbox { for all } l \in J \hbox{ and }  x \in
\om_{S/I}.
\] 
In the following we identify $\om_{S/I}$ with its image under $a$.

Let
\begin {equation}  \label{eq!KMlequ1}
                     {\bf L} \to S/I, \quad   
                     {\bf M}  \to S/J   
\end{equation} 
be minimal resolutions as $S$-modules. According to~\cite{FOV}
Proposition~A.2.12,
the dual complexes 
\begin{equation}    \label{eq!KMlequ2}
     {\bf L^*} \to \om_{S/I},  \quad {\bf M^*} \to \om_{S/J}
\end{equation}
are also minimal resolutions, where $*$  means $\Hom( \blank , S)$.
More precisely we have an exact sequence 
\begin{equation}   \label{eq!KMlequ3} 
     M_{r}^*  \xrightarrow{b}  M_{r+1}^* 
    \xrightarrow {c}  \om_{S/J} \to 0  
\end{equation}
and we set
\begin {equation}   \label{dfn!KMldefnofT}
               T = \ker c.
\end{equation}
For simplicity of notation, since $M_{r+1}^*$ is free of rank (say)  $l$,
 we identify it with $S^l$, so we have an exact  sequence 
\begin{equation}     \label{eq!KMlequ4}
    0 \to T \to S^l \xrightarrow{c}  \om_{S/J} \to 0.
\end{equation} 
For the canonical base $e_i=(0,\dots ,1,\dots ,0)\in S^l$, we fix liftings
\begin{equation}  \label{eq!KMlequ5} 
      q_i \in \Hom (J, \om_{S/I})    
\end{equation} 
of the basis  ${\overline e_i}=c(e_i) \in \om_{S/J}$, under the map $\res$
of
(\ref{eq!KMladjunct}).  
Notice that $\Hom (J, \om_{S/I})$ is generated by $\om_{S/I}$ together
with 
$q_1, \dots , q_l$. Moreover, clearly 
\begin{equation}  \label{eq!KMwhatisT}
  T = \big\{ (b_1, \dots ,b_l) \in S^l : \sum_i {b_iq_i} \in \om_{S/I}
\big\}.
\end{equation}

We denote by  $s_i \colon J \to T$ the map with 
\[
   s_i(t) = (0 , \dots , 0,t,0, \dots, 0),
\]
$t$ in the $i$th coordinate.  Then we define 
\[ 
      \Phi \colon  \Hom(T,\om_{S/I} ) \to \Hom(J,\om_{S/I})^l
\]
with 
\[
      \Phi (e) = (e\circ s_1, \dots , e\circ s_l).
\]

\begin{lemma}   \label{lemma!KMlinjimag}
The map $\Phi$ is injective, with  image equal to 
\[
    L  = \big\{ (k_1, \dots ,k_l) : \sum_i b_ik_i \in \om_{S/I} 
      \hbox{ whenever } (b_1, \dots ,b_l) \in T \big\}.
\]
\end{lemma}

\begin{pf} 
First of all notice that since $S/I$ is Cohen--Macaulay, there exists $t
\in J$
that is $S/I$-regular. Since $\om_{S/I}$ is a maximal Cohen--Macaulay
$S/I$-module,
$t$ is also $\om_{S/I}$-regular (compare e.g.~\cite{Ei} p.~529). 
Assume $e\circ s_i = e'\circ s_i$ for all $i$ and let  $b'= (b_1,\dots
,b_l) \in
T$. 
Then 
\[
   te(b')=e(tb')= \sum_i b_ie\circ s_i(t) = \sum_i b_ie'\circ s_i(t)=
te'(b'),
\]
and since $t$ is $\om_{S/I}$-regular we have $e=e'$. 
Moreover, this also shows that the image of $\phi$ is contained in $L$.

Now consider $(k_1,\dots , k_l) \in L$. Define $e \colon T \to \om_{S/I}$
with
$e(b_1,\dots ,b_l) = \sum b_ik_i$. Then $e\circ s_i(t)=tk_i =k_i(t)$,  so
$\Phi(e)=
(k_1, \dots ,k_l)$.
\QED \medskip
\end{pf}

Now we present a method, originally developed in \cite{KM}, and prove that
it
calculates
a set of generators for
\[
    \Hom (J, \om_{S/I} )/\om_{S/I},
\]
which was conjectured by Reid.

The natural map $\ S/I \to S/J$ induces
a map of complexes $ \psi \colon  {\bf L} \to {\bf M}$ and the dual map 
$ \psi^* \colon  {\bf M}^* \to {\bf L}^*$.  Using (\ref{eq!KMlequ2}), we
get
a commutative diagram with exact rows
\begin{equation*}
\begin{CD}
     M_{r-1}^* @>>>  M_r^*  @>b>> M_{r+1}^* @>c>> \om_{S/J}   \\
     @VVV  @VVV   \\
     L_{r-1}^* @>>> L_r^* @>>> \om_{S/I} @>>> 0 
\end{CD}
\end{equation*}
therefore, by the definition of $T$, we have an induced map 
$\psi^* \colon T \to \om_{S/I}$. Notice that this map is not canonical, 
but depends on the choice of $\psi$; we fix one such choice.  Set 
\[
   \Phi(\psi^*) = (k_1, \dots ,k_l). 
\]

\begin {theorem}   \label{thm!Kmlbasic}
The $S$-module  $\Hom_S(J, \om_{S/I})$ is generated by $\om_{S/I}$
together with $k_1, \dots, k_l$.
\end{theorem}

\begin {pf}  
Since $\om_{S/I}$ together with the $q_i$ of (\ref{eq!KMlequ5})
generate $\Hom(J, \om_{S/I})$, there are  equations
\begin {equation} \label{eq!KMlequ6}
  k_i  =  \sum_j a_{ij} q_j + \theta_i,  \quad a_{ij} \in S, \ \theta_i
\in
\om_{S/I}.
\end{equation}
Clearly $(q_i), (\theta_i)$ are in the image of $\Phi$, which by 
Lemma~\ref{lemma!KMlinjimag}
is equal to $L$; set  $(q_i) = \Phi(Q), \  (\theta_i) = \Phi(\Theta )$.
Define the map $f_0 \colon S^l \to S^l$ with 
\[
    f_0 (b_1, \dots b_l) = (b_1, \dots ,b_l) [a_{ij}].
\] 
 
Using (\ref{eq!KMwhatisT}), $T$ is invariant under $f_0$,
so there is an induced  map $f_1 \colon T \to T$ and, 
using (\ref{eq!KMlequ4}), a second induced map  $f_2 \in \Hom (\om_{S/J},
\om_{S/J})$,
with $f_2 ({\overline e_i}) = \sum_j a_{ij} {\overline e_j}$.  
We will  show that $f_2$ is an automorphism of $\om_{S/J}$, which will
prove 
the theorem.

Using (\ref{eq!KMlequ1}), (\ref{eq!KMlequ2}) and the definition 
(\ref{dfn!KMldefnofT}) of $T$, we get 
\begin{eqnarray*}
   \Ext^{r+1}(\om_{S/J},S) & = & \Ext^{r}(T,S) = S/J, \\ 
    \Ext^{r}(\om_{S/I},S) & = & S/I. 
\end{eqnarray*}
Moreover, by \cite{BH} Theorem~3.3.11 the natural map
$S/J \to \Hom (\om_{S/J}, \om_{S/J})$ is an isomorphism.  This, together
with 
the formal properties of the $\Ext$ functor,  imply that the natural map 
\[
  \Ext^{r+1}(\blank ,\id_S) \colon  \Hom (\om_{S/J}, \om_{S/J} ) \to S/J =
\Hom 
(S/J,S/J)
\]
is the identity $S/J \to S/J$.  
Therefore it is enough to show that  $\Ext^{r+1}(f_2)$ is a unit in $S/J$, 
and  by (\ref{eq!KMlequ4}) $\Ext^{r+1}(f_2) = \Ext^{r}(f_1)$, where 
for simplicity of notation we denote $\Ext^{*}(\blank ,\id_S)$ by
$\Ext^{*}(\blank )$.

By (\ref{eq!KMlequ6}) and the injectivity of $\Phi$ 
(Lemma~\ref{lemma!KMlinjimag})
\begin {equation} \label{eq!KMlequ7}
       \psi^* = Q\circ  f_1 + \Theta; 
\end{equation}
hence
\[
  \Ext^r(\psi^*) = \Ext^r(f_1) \Ext^r(Q) + \Ext^r (\Theta) 
\]
as maps $S/I \to S/J$.
Since $\Theta$ can be extended to a map $S^l \to \om_{S/I}$, 
$\Ext^r (\Theta) = 0$. Moreover, by the construction of $\psi^*$, 
$\Ext^r (\psi^*) = 1 \in S/J$. This  implies that  $\Ext^r (f_1)$
is a unit in $S/J$, which finishes the proof.  
\QED \medskip
\end{pf}

The arguments in the proof of Theorem~\ref{thm!Kmlbasic} also prove the
more 
general

\begin {theorem}  \label{thm!Kmlbasic2}
Let 
\[
f \colon T \to \om_{S/I}
\]
be an $S$-homomorphism, and set $f_i=f \circ s_i$ for  $1 \leq i \leq l$.
Then
$\om_{S/I}$ together with $f_1, \dots ,f_l$ generate $\Hom(J, \om_{S/I})$
if and only if the map
\[
    \Ext^r(f) : S/I \to S/J 
\]
is surjective.
\end{theorem}

\section {Unprojection of a complete intersection inside a complete
           intersection}   \label{sec!uciinsideci}

Let $S$ be a Gorenstein local ring  and $I \subset
J $ ideals of $S$, of codimensions $r$ and $r+1$ respectively.
We assume that  each is generated by a regular sequence, say
\begin {equation} \label{eq!forIJ}
  I = (v_1, \dots , v_r), \quad J = (w_1, \dots ,w_{r+1}).
\end{equation}
Since $I \subset J$, there exists an  $r \times (r+1)$ matrix $Q$ with
\begin {equation}  
   \begin {pmatrix} 
         v_1 \\ \vdots \\ v_r  
   \end {pmatrix}   = Q
     \begin {pmatrix} 
         w_1 \\ \vdots \\ w_{r+1}  
   \end {pmatrix}.   
\end {equation}

\begin{defn} \label{defn!ofwedgeofA}  
 $ \bigwedge^r Q $  is  the $1 \times (r+1)$ matrix
whose $i$th entry $(\bigwedge^r Q)_i$ is $(-1)^{i+1}$ times the
determinant 
of the submatrix of $Q$ obtained by removing the $i$th column.
\end{defn}

\begin{lemma} [Cramer's rule]  \label{lem!cramer}
For all $i,j$ the element  
\[
   (\bigwedge^r Q)_iw_j - (\bigwedge^r Q)_jw_i
\]  
belongs to the ideal $(v_1, \dots ,v_r)$.
\end{lemma}
\begin {pf} Simple linear algebra. \QED \medskip
\end{pf}

We define $g_i \in S$ by
\[
    \bigwedge^r Q = (g_1, \dots , g_{r+1}).
\]

The special case $r=2$ of the following theorem was proven by direct 
methods in \cite{CFHR} Lemma~6.11.

\begin {theorem}  \label{th!unpciinci}
$\Hom_{S/I}(J/I, S/I)$
is generated as $S/I$-module by two elements $\id$ and $s$, where
\[
   s(w_i) = g_i, \quad \hbox { for } 1 \leq i \leq r+1.
\]
\end{theorem}

\begin{pf}
Since $S/I$ is Gorenstein, we have $\om_{S/I}=S/I$.
By Lemma~\ref{lem!cramer} $s$ is well defined. Consider the minimal 
Koszul complexes corresponding to the generators given 
in (\ref{eq!forIJ}) that resolve $S/I,S/J$ as $S$-modules,
\[
                     {\bf L} \to S/I, \quad    {\bf M}  \to S/J. 
\]
The matrix $Q$ can be considered as a map $L_1 \to M_1$ making the
following
square commutative
\begin{equation*} 
\begin{CD}
                 L_1 @>(v_1, \dots, v_r)>> L_0=S  \\
                  @VQ^t VV   @V\id VV   \\
                 M_1  @>(w_1,\dots,w_{r+1}) >> M_0=S 
\end{CD}
\end{equation*}
 There are induced maps   $\bigwedge^n Q^t \colon L_n \to M_n$, 
$\bigwedge^0Q^t=\id_S$ (compare e.g.~\cite{BH} Proposition~1.6.8), 
giving a commutative diagram
\begin{equation*}
 \begin{CD}
                {\bf L} @>>> S/I  \\
                  @VVV   @VVV   \\
                {\bf M } @>>> S/J 
  \end{CD}
\end{equation*}
The last nonzero map is given by  $\bigwedge^r Q^t$; therefore using 
Theorem~\ref{thm!Kmlbasic}  the result follows. \QED 
\end{pf}

\begin {rem} Using \cite{PR} Theorem~1.5, the ideal $(v_i , Tw_j-g_j)$ of 
the polynomial ring  $S[T]$ is Gorenstein. It will be interesting
to find a geometric proof of this result. 
\end{rem}

\section { Tom \& Jerry }    \label{sec!tomandjerry}

\subsection {Definition of Tom}  \label{subs!tjdfntom}

Assume $S$ is a commutative Noetherian ring and  $a_{ij}^k,x_k,z_k  \in S$
for
$1 \leq k \leq 4,\ 2 \leq i  < j \leq 5$, and consider the $5 \times 5$
skewsymmetric matrix 
\begin{equation}   \label{eq!dfnoftommatr}
    A =  \begin {pmatrix}
             . & x_1 & x_2  & x_3  & x_4    \\
               &  .  & a_{23} &a_{24}  & a_{25}  \\
               & & . & a_{34} & a_{35} \\
               & \skewentry  &&.&a_{45}\\
               &&&&.
      \end{pmatrix}
\end{equation}
where
\[
     a_{ij} = \sum_{k=1}^4 a_{ij}^kz_k.
\]
The {\em Tom ideal} corresponding to $A$ is the ideal $I \subset S$
generated by
the 
Pfaffians of $A$. Notice that $I \subset J = (z_1, \dots ,z_4)$, as we
shall see
this is 
(under genericity conditions) an unprojection pair.

If $S$ is the polynomial ring over the ring of integers $\Z$ with
indeterminates 
$a_{ij}^k,x_k,z_k$ (indices as above) we call $I$ the {\em generic
integral Tom}
ideal.

\subsection {The generic integral Tom ideal is prine}

\begin {theorem} \label{thm!tomprimgor}
The generic integral Tom ideal $I$ is prime  of codimension three and
$S/I$ is
Gorenstein.
\end {theorem}

We prove the theorem  following the ideas of \cite{BV}~Chapter~2. 
Set $R=S/I$ and $X= \Spec R$. We will prove that $R$ is a domain.

\begin{lemma} \label{lem!localisdomain}
For all $i =1, \dots ,4$ the element $z_i \in R$ is not nilpotent.
Moreover, 
$R[z_i^{-1}]$ is a domain and $\dim R[z_i^{-1}] = 29$.
\end{lemma}

\begin{pf}  Due to symmetry, it is enough to prove it for $z_1$. 
By the form of the generators of $I$ it follows immediately that
$z_1 \in R$ is not nilpotent.
Consider the ring
\[
            T = k[z_1][z_1^{-1}][x_1,\dots ,x_4,z_2,z_3,z_4][a_{ij}^k],
\]
and the two skewsymmetric matrices $N_1$ and $N_2$ with
\[
     N_1 =     \begin {pmatrix}
             . & x_1 & x_2  & x_3  & x_4    \\
               &  .  & a_{23}^1 &a_{24}^1  & a_{25}^1  \\
               & & . & a_{34}^1 & a_{35}^1 \\
               & \skewentry &&.&a_{45}^1\\
               &&&&.
      \end{pmatrix}
\]
and $N_2 = A$, the generic Tom matrix defined in (\ref{eq!dfnoftommatr}).
Denote by $I_i$ the ideal of $T$ generated by the Pfaffians of $N_i$ for 
$i=1,2$. Consider the automorphism $f \colon T \to T$ that is the identity
on $k[z_1][z_1^{-1}][x_1,\dots ,x_4,z_2,z_3,z_4], \ f(a_{ij}^t) =
a_{ij}^t$ if 
$t \not= 1 $  and  $f(a_{ij}^1) = \sum_{k=1}^4 a_{ij}^kz_k$ ($f$ is 
automorphism since $z_1$ is invertible in $T$).  
Because $N_1$ is the generic skewsymmetric matrix, $I_1$ is
prime of codimension $3$ (see e.g.~\cite{KL}). Hence $I_2= f(I_1)$ is
also prime of codimension $3$, which proves the lemma.
\QED \medskip
\end{pf}

Set $U_i = \Spec R[z_i^{-1}] \subset X$, by Lemma~\ref{lem!localisdomain}
$U_i$ is irreducible.  Since the prime ideal $(x_k,a_{ij}^k) \subset R$
is in the intersection of the $U_i$, we have that $V = \cup U_i$ is 
also irreducible of dimension $29$.  In addition, $X \setminus V=
\Spec S/J$ has dimension $28$, where $J=(z_1,\dots ,z_4) \subset S$. 
Therefore, $X$ has dimension $29$, so $I$ has codimension three. Since $I$ 
is generated by Pfaffians, Theorem~\ref{thm!grade3compexact} implies 
that $R=S/I$ is Gorenstein. 
It follows that $I$ is unmixed, so $J$ is not contained in any associated
prime of $I$.  

Since $R[z_1^{-1}]$ is a domain, there is exactly one associated prime
ideal $P$ of $R$ such that $z_1 \notin P$. If $P$ is the single associated
prime ideal of $R$, then $z_1$ is a regular element of $R$ and $R$ is also
a
domain.
Suppose there is a second associated ideal $Q \not= P$. By what we have 
stated above and since $z_1 \in Q$, there is some $z_i \notin Q$. 
Since $R[z_i^{-1}]$ is a domain, it follows as before that $z_i \in P$.  
Now $PR[z_1^{-1}]= 0$, but the image of $z_i$ in $PR[z_1^{-1}]$ is
different 
from $0$ (otherwise $z_iz_1^t \in I$ which clearly doesn't happen), 
a contradiction which finishes the proof of Theorem~\ref{thm!tomprimgor}.

\subsection {Fundamental calculation for Tom} \label{subs!tjfcftom}

In this section $I \subset S$ is the generic integral Tom ideal 
defined in Subsection~\ref{subs!tjdfntom}.

Explicitly, $I = (P_0, \dots ,P_4)$ with
\begin{eqnarray}   \label{eq!genTom}
  P_0 & = & a_{{23}}a_{{45}}-a_{{24}}a_{{35}}+a_{{25}}a_{{34}} \\
  P_1 & = & x_{{2}}a_{{45}}-x_{{3}}a_{{35}}+x_{{4}}a_{{34}}  \notag   \\
  P_2 & = & x_{{1}}a_{{45}}-x_{{3}}a_{{25}}+x_{{4}}a_{{24}}   \notag \\
  P_3 & = & x_{{1}}a_{{35}}-x_{{2}}a_{{25}}+x_{{4}}a_{{23}}  \notag\\
  P_4 & = & x_{{1}}a_{{34}}-x_{{2}}a_{{24}}+x_{{3}}a_{{23}}  \notag
\end{eqnarray}
Clearly $ I \subset J = (z_1,z_2,z_3,z_4)$.  

Since each $P_i$, for $1 \leq i \leq 4 $, is linear in $z_j$, there exists 
(unique) $4 \times 4 $ matrix $Q$ independent of the $z_j$ such that
\begin {equation}  \label{eq!tomlinear}
   \begin {pmatrix} 
         P_1 \\ \vdots \\ P_4  
   \end {pmatrix}   = Q 
     \begin {pmatrix} 
         z_1 \\ \vdots \\ z_4  
   \end {pmatrix}
\end {equation} 
We denote by  $Q_i$ the $i$th row of $Q$, and by  $\widehat{Q}_i$ the
submatrix of $Q$ obtained by deleting the $i$th row.  
Since (compare~(\ref{eq!pfafcompl}))
\[
       x_4P_4 = x_1P_1-x_2P_2+x_3P_3,
\] 
and $Q$ is independent of the $z_j$, it follows that
\begin{equation} \label{eq!forQ}             
       x_4Q_4 = x_1Q_1-x_2Q_2+x_3Q_3.
\end{equation} 

For $i=1,\dots ,4$ we define a $1 \times 4$ matrix $H_i$  by 
\[
    H_i =  \bigwedge^3 \widehat{Q}_i,
\]
where $\bigwedge$ as in Definition~\ref{defn!ofwedgeofA}.

\begin {lemma}  \label{lem!tomfirstid} 
For all $i,j$ 
\[
   x_i  H_j =  x_j H_i.
\]
\end {lemma}

\begin{pf} Equation (\ref{eq!forQ}) implies, for example,  
that  
\[  
  \bigwedge^3 \widehat{Q}_3  = 
   \bigwedge^3
   \begin {pmatrix}
             Q_1 \\ Q_2  \\ Q_4 
   \end{pmatrix}  = 
   \bigwedge^3
   \begin {pmatrix}
             Q_1 \\ Q_2  \\ \frac {x_3}{x_4} Q_3 
    \end{pmatrix} = 
     \frac {x_3}{x_4}   \bigwedge^3 \widehat{Q}_4.  
\]  \QED \medskip
\end{pf}

Using the previous lemma we can define four polynomials $g_i$ by
\begin{equation}    \label{eq!dfnoftomgi}
  (g_1, g_2, g_3, g_4)  =  \frac {H_j} {x_j} 
\end{equation} 
and this definition is independent of the choice of $j$. 

\begin {lemma}  \label{lem!tomsecondid} 
For all $i,j$ 
\[
   g_i z_j - g_j z_i  \in I.
\]
\end {lemma}

\begin{pf} The definition (\ref{eq!tomlinear}) of $Q$ 
implies
\begin {equation*}  
   \begin {pmatrix} 
         P_1 \\ P_2 \\ P_3   
   \end {pmatrix}   =   
   \begin {pmatrix} 
         Q_1 \\ Q_2 \\ Q_3   
   \end {pmatrix} 
     \begin {pmatrix} 
         z_1 \\ \vdots \\ z_4  
   \end {pmatrix}   
\end {equation*} 
By Cramer's rule (Lemma~\ref{lem!cramer})
\[
    (H_4)_i z_j -  (H_4)_j z_i \in I,
\]  
as a consequence
\[
    x_4 (g_i z_j -g_j z_i ) \in I.
\] 
$I$ is prime by Theorem~\ref{thm!tomprimgor}, so the result follows.
\QED \medskip
\end{pf}

\subsection {Resolution of Tom}    \label{subs!tjcomplexes}

We use the notations of Subsection~\ref{subs!tjfcftom}.

Consider the Koszul complex {\bf M} that gives a resolution of
the ring $S/J$
\[
   0 \to S \xrightarrow{B_4} S^4 \xrightarrow{B_3} S^6  
    \xrightarrow{B_2} S^4 \xrightarrow{B_1} S \to 0,      
\]
with  
\begin{eqnarray*}
   B_1  &=&  (z_1, z_2, z_3, z_4),   \\
   B_2  &=&  \begin{pmatrix}  
   -z_2 &  -z_3 & 0 &  -z_4 & 0 & 0 \\
   z_1 &  0 & -z_3 &  0 &  -z_4 &  0 \\
   0 &  z_1 & z_2 & 0 & 0 & -z_4 \\
   0 & 0 & 0 & z_1 & z_2 & z_3
      \end{pmatrix}
\end{eqnarray*}
and
\begin{eqnarray*}
   B_3 =  \begin{pmatrix} 
     z_3 & z_4 & 0 & 0 \\ 
    -z_2& 0& z_4& 0 \\ 
    z_1& 0& 0& z_4 \\ 
    0& -z_2& -z_3& 0 \\ 
    0& z_1& 0& -z_3 \\ 
    0& 0& z_1& z_2 
         \end{pmatrix}, \
   B_4 = (-z_4, z_3, -z_2, z_1) ^t.
\end{eqnarray*}

Moreover, the skewsymetric matrix $A$ defines as 
in (\ref{eq!pfafcompl}) a complex {\bf L}:
\begin{equation}   \label{eq!tomfdresofI} 
  0 \to S \xrightarrow{C_3} S^5 \xrightarrow{C_2} S^5 \xrightarrow{C_1} S
\to 0      
\end{equation}
resolving the ring $S/I$. 
Here $C_2 = A,  C_1 = (P_0, -P_1, P_2, -P_3, P_4)$,
and $C_3$ is the transpose matrix of $C_1$. 
Define the $4 \times 1$ matrix  $D_3$ with 
\[
   D_3 = (-g_4, g_3, -g_2, g_1)^t, 
\]
where the $g_i$ are as in (\ref{eq!dfnoftomgi}). 

\begin{theorem}   \label{thm!cmltomcd}
There exist matrices $D_2,D_1,D_0$ (of suitable sizes) making the 
following diagram commutative.
\begin{equation*}
\begin{CD}
     0  @>>>  S   @>C_3>> S^5  @>C_2>> S^5 @>C_1>> S   \\
      @. @VD_3VV  @VD_2VV  @VD_1VV @VD_0VV  \\
     S @>B_4 >> S^4  @>B_3>> S^6 @>B_2>> S^4 @>B_1>> S  
\end{CD}
\end{equation*}
In addition we can assume that $D_0 = 1 \in \Z$.  
\end{theorem}

\begin{pf} 
As in (\ref{eq!KMlequ2}), the dual complexes 
\[
   S^* \to  S^{5*} \to S^{5*} \to S^* 
\] 
and
\[
   S^* \to S^{4*} \to S^{6*} \to S^{4*} 
\]
are exact. Using  Lemma~\ref{lem!tomsecondid}, there exists
$D_2^*$ making the (dual) square commutative. Then, the existence 
of $D_1^*$ and $D_0^*$ follows by simple homological algebra.

We get $D_0 \in \Z$  by checking the degrees in 
the commutative diagram.  The fact that we can take $D_0 = 1$ will
be proved in the next subsection.
\QED \medskip
\end{pf}

\subsection {Original Tom} 

Write $Y \subset \Pp^8$ for the Segre embedding of  $\Pp^2 \times \Pp^2$
(all schemes are over the ring of integers $\Z$).
It is easy to see that $Y$ is projectively Gorenstein of codimension
four.  The defining  equations are $\rank L \leq 1$, where $L$ is the
generic $3 \times 3$ matrix
\[ 
     L =   \begin {pmatrix}
         a & x_3 &  x_4  \\
          x_1 & z_1 & z_2 \\
         x_2 & z_3 &  z_4  
    \end {pmatrix}
\]
and $a,x_1, \dots ,x_4, z_1, \dots z_4$ are indeterminates. Let
$X \subset \Pp^7$ be the image of the projection of $Y$ from the point 
$((1,0,0),(1,0,0)) \in \Pp^2 \times \Pp^2$. Clearly, the ideal of $X$ is 
generated by the five polynomials (the five minors of $L$ {\it not} 
involving $a$)
\[
    I(X)=(x_3z_2-x_4z_1,x_3z_4-x_4z_3,x_1z_3-x_2z_1,   
       x_1 z_4- x_2 z_2).
\] 
These are the Pfaffians of the skewsymmetric $5 \times 5$ matrix
\begin {equation}  \label{eq!chp3prortom}
     M = \begin {pmatrix}
           . & x_1 & x_2 &  x_3 & x_4 \\
             & . & 0 & z_1 & z_2 \\
             &   &  . & z_3 & z_4 \\
             &   &   &  & 0 \\
             &   &   &    & .
         \end {pmatrix}
\end {equation}
We call $X$ the {\em original Tom}. $X$ contains the complete intersection 
$D$ (the exceptional locus of the inverse map to the projection 
$Y \broken X$), with $I(D)=(z_1, \dots ,z_4)$. 

Equations (\ref{eq!dfnoftomgi}) specialize to 
\[
   g_1 = x_1x_3, \ \ g_2 = x_1x_4, \ \ g_3 = x_2x_3, \ \ g_4 = x_2 x_4.
\]

An easy calculation using (the specialization of) the complexes {\bf L},
{\bf M}
defined in Section~\ref{subs!tjcomplexes} gives that we can take in the 
specialization of the diagram of Theorem~\ref{thm!cmltomcd}
\[
  D_2' =  \begin{pmatrix} 
        0 &-x_2 & 0 & 0 & 0 \\
        0 & 0 & 0 & -x_4 & 0 \\
        0 & 0 & x_2 & x_3 & 0 \\
        0 & x_1 & 0 & x_3 & 0  \\ 
        0 & 0 & 0 & 0 & x_3 \\
        0 & 0 & x_1 & 0 & 0
          \end{pmatrix},
\]
\[
   D_1' = \begin{pmatrix}
    -z_4 & 0 & x_4 & 0 &  -x_2 \\ 
    0&0&-x_3&x_2&0 \\
   z_2&-x_4&0&0&x_1\\ 
   0&x_3&0&-x_1&0 
      \end{pmatrix}
\]    
and  $D_0' =1$.  

\paragraph {End of proof of Theorem~\ref{thm!cmltomcd}}
Using the uniqueness up to homotopy of a map between resolutions of
modules 
induced by a fixed map between the modules, the last part of 
Theorem~\ref{thm!cmltomcd} follows from $D_0'=1$.

\subsection {Local Tom}  
 
Let $S$ be a Gorenstein local ring, $a_{ij}^k  \in S$ and $x_k,z_k \in m$,
the
maximal ideal of $S$, with  indices as in  Subsection~\ref{subs!tjdfntom}. 
Let $A$ be the skewsymmetric matrix 
defined in (\ref{eq!dfnoftommatr}), $I$ the ideal generated by
the  Pfaffians  of $A$ (see~(\ref{eq!genTom})) and $J=(z_1, \dots ,z_4)$.

We assume that $z_1, \dots ,z_4$ is a 
regular sequence and that $I$ has codimension three, the maximal possible.
Since $S$ is Cohen--Macaulay, the grade of $I$ is also three.
By Theorem~\ref{thm!grade3compexact},  the complex 
${\bf L}$ defined in (\ref{eq!tomfdresofI}) is the minimal resolution
of $S/I$ and $S/I$ is Gorenstein.

Recall that in (\ref{eq!dfnoftomgi})  we defined elements $g_i$ which are 
polynomials of  $a_{ij}^k$ and $x_k$. Define a map  $\psi \colon J/I \to
S/I$
with $z_i \mapsto g_i$. By $\res$ we denote the residue 
map defined in (\ref{eq!KMladjunct}).

\begin {theorem}  \label{thm!localtomthm}
The element $\res (\psi) \in S/J$ is a unit, and the ideal 
\begin {equation}  \label{eq!unprojtomideal}
   (P_0, \dots ,P_4, Tz_1-g_1, \dots , Tz_4-g_4)
\end {equation} 
of the polynomial ring $S[T]$ is Gorenstein of codimension four.
\end{theorem}

\begin{pf}
The theorem follows immediately from Theorem~\ref{thm!cmltomcd} (since the
diagram is defined over $\Z$), Theorem~\ref{thm!Kmlbasic}, and 
\cite{PR}~Theorem~1.5.
\QED \medskip
\end{pf}

\begin{rem} It is easy to see that in the case of the unprojection of the
generic
integral Tom, or more generaly if $S$ is arbitrary and $I$ is sufficiently
general, 
the unprojection ideal $(P_i, Tz_i-g_i) \subset S[T]$ is generated by 
$9$ elements and has $16$ syzygies.
\end{rem}  

\begin{exa} We will calculate the equations of the unprojection of
the Tom matrix 
\[
   A =    \begin {pmatrix}
             . & x_1 & x_1x_2z_2  & x_3   & x_3    \\
               &  .  &  0 & x_3z_3z_1+z_2  & z_3   \\
               & & . &  z_4 &  z_1 \\
               &\skewentry &&.& 0 \\
               &&&&.
          \end{pmatrix}
\]
We make the choice $a_{24}^1 = x_3z_3, a_{24}^3 =0$ and replace $A$ 
with
\[
   A' =   \begin {pmatrix}
             . & t_1 & t_2  & t_3  & t_4    \\
               &  .  &  0 & x_3z_3z_1+z_2  & z_3   \\
               & & . &  z_4 &  z_1 \\
               &\skewentry &&.& 0 \\
               &&&&.
           \end{pmatrix}
\]
We have 
\[
      Q' =  \begin{pmatrix}
              -t_3 & 0 & 0 & t_4 \\
             t_4x_3z_3 & t_4 & -t_3 & 0 \\
               t_1 & 0 & -t_2 & 0 \\
              -t_2x_3z_3 & -t_2 & 0 & t_1 
             \end{pmatrix}
\]
and 
\[
 H_4' =( t_4t_4t_2,\  t_4(-t_4x_3z_3t_2+t_1t_3),\ t_4t_1t_4,\  t_4t_2t_3).  
\] 
Therefore, the unprojection equations are
\begin{eqnarray*}
  \begin {matrix} 
    g_1  = x_3(x_1x_2z_2), &  g_2  =  -x_3x_3z_3(x_1x_2z_2)+x_1x_3, & 
        g_3 = x_1x_3, \\  
    g_4 = (x_1x_2z_2)x_3. & & 
  \end{matrix}
\end{eqnarray*}
\end{exa}

\begin {exa} (Del Pezzo surfaces)
Using serial unprojection, we can explicitly construct
the equations of del Pezzo surfaces. For example
start with the smooth cubic hypersurface
$X_3 \subset \P^3$ with equation $xz(x+z)-yw(y+w)=0$
over $\C$. It is a general fact (see e.g. \cite{Har}
Section~V.4)
that it contains six lines which are mutually disjoint, 
for our purposes we write down four lines which
are mutually disjoint:
\[
  \begin {matrix}
     l_1: x=y=0 &   l_2: z = w =0 \\
     l_3: x+z = y+w =0  & l_4: x-w=x+y+z=0
  \end{matrix}
\]
Since $X_3$ is smooth, unprojection corresponds to
blowing down a $-1$ curve.

We unproject the line $l_1$ in $X_3$, to get the
del Pezzo degree $4$ surface $X_4 \subset \P^4$,
with equations $sx-w(y+w) = sy-z(x+z) =0$, where
$s$ is the new coordinate.

Then we unproject
the strict transform $l_2'$ of $l_2$ in
$X_4$. The line $l_2'$ has equations $s=z=w=0$.
Using Theorem~\ref{th!unpciinci}, the unprojection  
is the degree $5$ del Pezzo $X_5 \subset \P^5$
with equations given by the five Pfaffians
$g_1, \dots , g_5$ of the skewsymmetric
$5 \times 5$ matrix
\[
  M = \begin{pmatrix}
        . & -x & w & -y & -z \\
          & . & 0     & t & -(y+w) \\
          &   & . & x+z & s \\
         & \skewentry &   &  .  & 0  \\                 
          &  &    &    &  .
      \end{pmatrix}
\]
 where $t$ is  the new coordinate.
This is a Tom format with respect to
the line $l_3' : t =s=   x+z= y+w= 0$
(the strict transform of $l_3$), so
we can unproject $l_3'$. Using
Theorem~\ref{thm!localtomthm}, we get
the del Pezzo surface $X_6 \subset \P^6$           
(a hyperplane section of the Segre embedding
$\P^2 \times \P^2 \subseteq \P^8$),
with ideal
\[
   I(X_6) = (g_1, \dots , g_5, ut-xy,
  u(y+w)+xz, u(x+z)+yw, us+wz),
\]
where $u$ is the new coordinate.
We leave it as an exercise to the interested
reader to find the equations of the
del Pezzo surface  $X_7 \subset \P^7$ 
arising as  unprojection   
of the strict transform of $l_4$ inside $X_6$
(compare also  Remark~\ref{thm!localtomthm}).
\end{exa}

\subsection {Jerry}
The Jerry unprojection family has many similarities with the Tom family,
we will
avoid
repeating the arguments when they are almost identical.

Assume $S$ is a commutative Noetherian 
ring and  $x_i,z_k,c^k,a_i^k,b_i^k \in S$ for
$1 \leq  k \leq 4, \ 1 \leq i \leq 3$, and consider the $5 \times 5$ 
skewsymmetric matrix 
\begin{equation}   \label{eq!dfnofjerrymatr}
    A =     \begin {pmatrix}
             . & c & a_1  & a_2   & a_3    \\
               &  .  & b_1 &  b_2  & b_3  \\
               & & . & x_1 & x_2 \\
               &\skewentry &&.& x_3 \\
               &&&&.
      \end{pmatrix}
\end{equation}
where 
\[
   a_i = \sum_{k=1}^4 a_i^kz_k, \ \    b_i =  \sum_{k=1}^4 b_i^kz_k, \ \
   c =  \sum_{k=1}^4 c^kz_k.
\]

The {\em Jerry} ideal corresponding to $A$ is the ideal $I \subset S$
generated by the Pfaffians of $A$. Notice that $I \subset J = (z_1, \dots
,z_4)$.

If $S$ is the polynomial ring over the ring of integers $\Z$ with
indeterminates 
$x_i,z_k,c^k,a_i^k,b_i^k$ (indices as above) we call $I$ the {\em generic
integral Jerry} ideal.

The proof of the following theorem is very similar to the proof of
Theorem~\ref{thm!tomprimgor}.
\begin {theorem} \label{thm!jerryprimgor}
The generic integral Jerry ideal $I$ is prime  of codimension three and
$S/I$ is
Gorenstein.
\end {theorem}

We work with the generic integral Jerry ideal $I$.  
Write $I = (P_1, \dots ,P_5)$, with 
\begin{eqnarray}   \label{eq!genJerry}
   P_1 &=&  b_1x_3-b_2x_2+b_3x_1  \\
   P_2 &=&  a_1x_3-a_2x_2+a_3x_1 \notag \\
   P_3 &=&  cx_3-a_2b_3+a_3b_2  \notag \\
   P_4 &=&  cx_2-a_1b_3+a_3b_1  \notag \\
   P_5 &=&  cx_1-a_1b_2+a_2b_1  \notag 
\end{eqnarray}

Unlike in the Tom case, we only have two Pfaffians, $P_1$ and  $P_2$,
linear
in $z_k$. $P_3$ is quadratic in $z_k$, but after choosing to consider
$a_2, a_3$ as indeterminates it can be considered linear.  Using this
convention we write
\begin {equation}  \label{eq!jerlinear}
   \begin {pmatrix} 
         P_1 \\ P_2 \\ P_3  
   \end {pmatrix}   = Q 
     \begin {pmatrix} 
         z_1 \\ \vdots \\ z_4  
   \end {pmatrix}
\end {equation}  
$Q$ is a $3 \times 4$ matrix,  with 
\begin{eqnarray*} 
    Q_{1k} & = &  b_1^k x_3 - b_2^k x_2 + b_3^k x_1  \\
    Q_{2k} & = &  a_1^k x_3 - a_2^k x_2 + a_3^k x_1  \\
    Q_{3k} & = &  c^k x_3 - a_2 b_3^k  +  a_3 b_2^k
\end {eqnarray*}   
We define $h_i$ by
\[ 
   \bigwedge^3 Q = (h_1, \dots ,h_4 ).
\]  
($\bigwedge$ as in  Definition~\ref{defn!ofwedgeofA}.)

\begin{lemma}  \label {lem!jerryfund}
For  $i=1,\dots ,4$ there  exist polynomials $K_i, L_i$ with
\[
      h_i =  x_3 K_i+ (a_2x_2-a_3x_1) L_i.
\] 
Therefore, we can write 
\[
    h_i = x_3 (K_i+a_1L_i) - L_iP_2.
\]
\end{lemma}

\begin{pf} Let $M$ be the matrix obtained from $Q$ 
by substituting $x_3=0$. Since 
\begin {eqnarray*} 
    M_{1k} & = &  x_1 b_3^k   - x_2 b_2^k \\
    M_{3k} & = &  -a_2 b_3^k + a_3 b_2^k 
\end{eqnarray*}  
we get
\[ 
     M =   \begin {pmatrix}
         x_1 & 0 &  -x_2  \\
          0 & 1 & 0 \\
         -a_2 & 0 &  a_3 
    \end {pmatrix} 
          \begin {pmatrix}
          b_3^1 & b_3^2 & b_3^3 & b_3^4  \\
           M_{21} & M_{22}  & M_{23} & M_{24} \\
          b_2^1 & b_2^2 & b_2^3 & b_2^4 
    \end {pmatrix}       
\]
The lemma follows from elementary properties of determinants. 
\QED \medskip
\end{pf}
   
We fix the polynomials $K_i, L_i$ defined (implicitly) in the proof of \\
Lemma~\ref{lem!jerryfund}. For $i=1, \dots ,4$  we define 
polynomials $g_i$ by
\begin {equation} \label{eq!dfnofjerrygi}
    g_i= K_i+a_1L_i.
\end{equation} 

Using (\ref{eq!jerlinear}) and Theorem~\ref{thm!jerryprimgor},
the following lemma follows in the same way as
Lemma~\ref{lem!tomsecondid}.

\begin {lemma}  \label{lem!jerrysecondid} 
For all $i,j$ we have $g_i z_j - g_j z_i  \in I$.
\end {lemma}

Consider as in Subsection~\ref{subs!tjcomplexes}
the Koszul complex {\bf M} resolving $S/J$, the complex {\bf L} resolving
the ring $S/I$, and define the $4 \times 1$ matrix  $D_3$ with 
\[
   D_3 = (-g_4, g_3, -g_2, g_1)^t, 
\]
where the $g_i$ are as in (\ref{eq!dfnofjerrygi}). 

Except from the part $D_0=1 \in \Z$, the following theorem follows
immediately
using the arguments in the proof of Theorem~\ref{thm!cmltomcd}.

\begin{theorem}   \label{thm!cmljerrycd}
There exist matrices $D_2,D_1,D_0$ (of suitable sizes) making the 
following diagram commutative.
\begin{equation*}
\begin{CD}
     0  @>>>  S   @>C_3>> S^5  @>C_2>> S^5 @>C_1>> S   \\
      @. @VD_3VV  @VD_2VV  @VD_1VV @VD_0VV  \\
     S @>B_4 >> S^4  @>B_3>> S^6 @>B_2>> S^4 @>B_1>> S  
\end{CD}
\end{equation*}
Moreover, we can assume that $D_0 = 1 \in \Z$.  
\end{theorem}     

The part $D_0=1$ follows,  as in the Tom case,  by a specialization
argument
using the 
original Jerry which we now define.

\paragraph{Original Jerry} 
Write $Y \subset \Pp^7$ for the image of the Segre embedding of  $\Pp^1
\times
\Pp^1
\times \Pp^1$. $Y$ is projectively Gorenstein of codimension
four.  Let $X \subset \Pp^6$ be  the image of the projection of $Y$ from
the
point 
$((1,0),(1,0),(1,0)) \in \Pp^1 \times \Pp^1 \times \Pp^1$.
If $z_1,\dots ,z_4, x_1, \dots ,x_3$ are homogeneous coordinates 
for $\P^6$, the homogeneous ideal of $X$ is given (compare \cite{Ki},
Example~6.10)
by the Pfaffians of the skewsymmetric matrix
 \[
     M = \begin {pmatrix}
           . & z_1 & z_2 &  z_3 & 0   \\
             & .   & 0   &  z_3 & z_4 \\
             &     &  .  &  x_1 & x_2 \\
             &     &     &      & x_3 \\
             &     &     &      & .
         \end {pmatrix}
\]
We call $X$ the {\em original Jerry}.

\paragraph{Local Jerry}

Let $S$ be a Gorenstein local ring, $a_i^k,b_i^k,c^k  \in S$ and $x_i,z_k
\in
m$,
the maximal ideal of $S$,  with indices as above. 
Let $A$ be the skewsymmetric matrix 
defined in (\ref{eq!dfnofjerrymatr}), $I$ the ideal generated by
the  Pfaffians  of $A$ (see~(\ref{eq!genJerry})) and $J=(z_1, \dots
,z_4)$.

We assume that $z_1, \dots ,z_4$ is a 
regular sequence and that $I$ has codimension three, the maximal possible. 
Since $S$ is Cohen--Macaulay, the grade of $I$ is also three.
By Theorem~\ref{thm!grade3compexact},  the complex 
${\bf L}$ defined in (\ref{eq!pfafcompl}) is the minimal resolution
of $S/I$ and $S/I$ is Gorenstein. 

Recall that in (\ref{eq!dfnofjerrygi})  we defined elements $g_i$ which
are 
polynomials in  $a_i^k,b_i^k,c_i^k,z_k$ and $x_i$. Define a map  
$\psi \colon J/I \to S/I$ with $z_i \mapsto g_i$. By $\res$ we denote the
residue 
map defined in (\ref{eq!KMladjunct}).

The proof of the following theorem is very similar to the proof of 
Theorem~\ref{thm!localtomthm}.

\begin {theorem}  \label{thm!localjerrythm}
The element $\res (\psi) \in S/J$ is a unit, and the ideal 
\[
   (P_1,\dots ,P_5, Tz_1-g_1, \dots, Tz_4-g_4)
\]
of the polynomial ring $S[T]$ is Gorenstein of codimension four.
\end{theorem}

\begin{rem} As in the Tom case, it is easy to see that in the case of the
unprojection of the generic
integral Jerry, or more generaly if $S$ is arbitrary and $I$ is
sufficiently
general, 
the unprojection ideal $(P_i, Tz_j-g_j) \subset S[T]$ is generated by $9$
elements
and has $16$ syzygies. As already mentioned in the introduction, it would
be
very interesting to
find out whether all known Gorenstein codimension four rings with $9
\times 16$
resolution
are related to the Tom and Jerry unprojection families. That would
probably need
an intrinsic
characterization of Tom and Jerry, which at present is not known.    
\end{rem}  

\section {Final comments on calculating unprojection}  \label{sec!fc}  

\begin{rem}
In all three families of unprojection data $ I \subset J \subset S$ 
we handled, $J=(w_1, \dots ,w_{r+1})$ was generated by a regular sequence, 
and the treatment was roughly as follows: Choose $r$ suitable elements 
$v_1, \dots , v_r$ of $I$, use Cramer rule with respect to the matrix 
$Q$ that expresses the $v_i$ as combination
of the $w_j$, and then (in the Tom and Jerry cases) divide by a common
factor. 
I believe that whenever $J$ is generated by a regular sequence the steps
will 
be similar, the question is to give a unified treatment of all such cases.
\end{rem}

\begin{rem}
The case where $J$ is not a complete intersection seems to
be more mysterious. Cramer's rule definitely does not work, as it can be
seen 
from \cite{KM}~Example~2.1.
\end{rem}

\begin{rem}
If one could give an interpretation of the residue epimorphism 
(\ref{eq!KMladjunct})
\[  
   \Hom_{S/I} (J/I, \om_{S/I}) 
             \xrightarrow{\, \res \,} \om_{S/J} \to 0
\] 
using differentials, then, in principle, he should be able to have a clean
and 
easy calculation of unprojection avoiding complexes.    
\end{rem} 

\begin{rem} 
All three families of unprojection we calculated belong to the type I
unprojection (see \cite{Ki}, Section~9), since the ring $S/J$ is
Gorenstein. Reid and his coworkers have discovered (at least) three more
types of unprojection, which they call Types II-IV, see loc.~cit.~and 
\cite{R2}. 
Since in Theorem~\ref{thm!Kmlbasic} we need $S/J$ to be Cohen--Macaulay, 
our methods can also  handle Type III unprojections, but need
improvement in order to calculate Type II and IV.    
\end{rem}

{\footnotesize \paragraph   {\it Acknowledgements} I am grateful to Miles
Reid
for 
suggesting the problem and many useful conversations, and to David Mond
and
Frank--Olaf Shreyer for important suggestions. This work is part of a
Warwick
PhD thesis, 
financially supported by the Greek State Scholarships Foundation.}

\bigskip
\noindent
Stavros Papadakis \\
Karaiskaki 11, Agia Paraskevi \\
GR 153 41, Attiki \\
Greece\\ 
email: stavrospapadakis@hotmail.com\\ 
web: www.maths.warwick.ac.uk/$\sim$miles/doctors/Stavros/

\end{document}